\documentclass[11pt]{article}
\usepackage{bm} 
\usepackage{graphicx}
\usepackage{epstopdf}
\usepackage{float}
\usepackage{amsthm}
\usepackage{natbib}
\usepackage{amsmath}
\newtheorem{theorem}{Theorem}

\makeatletter
\begin{document}
\title{Wasserstein Statistics in 1D Location-Scale Model}
\author{Shun-ichi Amari \\
RIKEN Center for Brain Science
}
\date{}
\maketitle

\begin{abstract}
Wasserstein geometry and information geometry are two important structures introduced in a manifold of probability distributions. The former is defined by using the transportation cost between two distributions, so it reflects the metric structure of the base manifold on which distributions are defined. Information geometry is constructed based on the invariance criterion that the geometry is invariant under reversible transformations of the base space.  Both have their own merits for applications. Statistical inference is constructed on information geometry, where the Fisher metric plays a fundamental role, whereas Wasserstein geometry is useful for applications to computer vision and AI. We propose statistical inference based on the Wasserstein geometry in the case that the base space is 1-dimensional. By using the location-scale model, we derive the $W$-estimator explicitly and studies its asymptotic behaviors.
\end{abstract}

\section{Introduction}

Wasserstein geometry defines a divergence between two probability distributions $p(x)$ and $q(x)$, $x \in X$ by using the cost of transportation from $p$ to $q$.  Hence, it reflects the metric structure of the underlying manifold $X$ on which probability distributions are defined.  Information geometry, on the hand, studies an invariant structure such that geometry does not change under transformations of $X$ which would change the distance within $X$.  So it is independent of the metric of $X$.

Both geometries have their own histories \citep*[see e.g.,][]{Villani2003,Villani2009,Amari2016}.  Information geometry has been successful for elucidating statistical inference, where the Fisher information metric plays a fundamental role.  It has successfully been applied, not only to statistics, but also to machine learning, signal processing, systems theory, physics and many others \citep*{Amari2016}.  Wasserstein geometry has been a useful tool for geometry, where the Ricci flow has played an important role \citep*{Villani2009, LM2018}.  Recently, it has a wide scope of applications in computer vision, deep learning and more \citep[e.g.,][]{FZMAP2015, ACB2917, MMC2015, PC2019}.  There are some trials to connect the two geometries.  \citet*{LZ2019} gave a unified theory connecting them.  See also \citet*{WL2019} and  \citet{AKO2018, AKOC2019}. 

It is natural to consider statistical inference from the Wasserstein geometry point of view and compare the results with information-geometrical inference based on the likelihood \citep*{LZ2019}.  The present short article studies the statistical inference based on the Wasserstein geometry from a different point of view of \citet*{LZ2019}.  Given a number of independent observations from a probability distribution belonging to a statistical model with a finite number of parameters, we define the $W$-estimator that minimizes $W$-divergence from the empirical distribution $\hat{p}(x)$ derived from observed data to the statistical model.  In contrast, the information geometry estimator is the one that minimizes Kullback-Leibler divergence from the empirical distribution to the model, and it is the maximum likelihood estimator.

We use 1D base space $X={\bm{R}}^1$, and define the transportation cost to be equal to the square of the Euclidean distance between two points in ${\bm{R}}^1$.  We further focus on the location-scale model to obtain explicit solutions in the asymptotic resume, that is, the number of observations is sufficiently large. We then give an explicit expression of the $W$-estimator, proving that it is asymptotically consistent and further calculate its asymptotic variance. Although they are not Fisher efficient, it minimizes the divergence between the empirical distribution and the model. We may say that it is $W$-efficient estimator in this sense.

The present $W$-estimator is different from \citet*{LZ2019}, based on the Wasserstein score function.  The $W$-efficiency of this estimator is defined.  Although this is a fundamental theory, opening a new paradigm connecting information geometry and $W$ geometry, it does not minimizes the $W$-divergence from the empirical one to the model.  It is an interesting problem to compare these two frameworks of Wasserstein statistics.

The present paper is organized as follows. After introduction, we formulate estimating equations for a general parametric statistical model in the 1D-case.  We show in section 2 that the optimal estimator uses only a linear function of observations.  We then focus on the location-scale model in section 3.  We give an explicit form of the $W$-estimator.  We analyze the asymptotic properties of the $W$-estimator. We studies the geometry of the location-scale model in section 4, showing that it is Euclidean \citep*{LZ2019}, although it is a curved submanifold in the function space of $W$-geometry \citep{Takatsu2011}.  We finally give characteristic features of the $W$-estimator, comparing it with the maximum likelihood estimator.

\section{$W$-estimator}

We first show the optimal transportation cost sending $p(x)$ to $q(x)$, $x \in {\bm{R}}^1$ when the transportation cost from $x$ to $y$, $x, y \in {\bm{R}}^1$, is $(x-y)^2$.  Let $P(x)$ and $Q(x)$ be the cumulative distributions of $p$ and $q$, respectively,
\begin{eqnarray}
  P(x) &=& \int^x_{-\infty} p(u) du, \\
  Q(x) &=& \int^x_{-\infty} q(u) du.
\end{eqnarray}
Then, it is known that the optimal transportation plan is to send mass of $p(x)$ at $x$ to $x'$, such that
\begin{equation}
  P^{-1}(x) = Q^{-1} \left(x' \right),
\end{equation}
$P^{-1}$ and $Q^{-1}$ being the inverse functions of $P$ and $Q$.
\begin{figure}
    \centering
    \includegraphics[width=8cm]{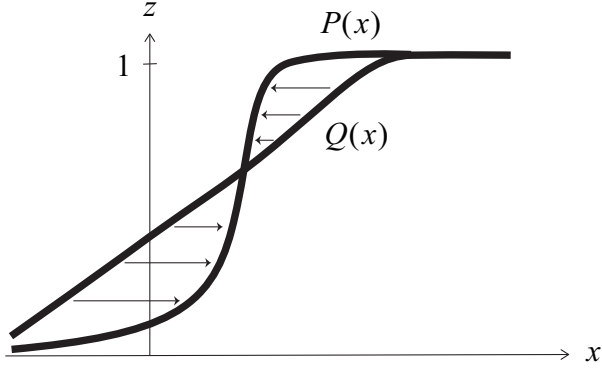}
    \caption{Optimal transportation plan from $p$ to $q$}
    \label{fig1}
\end{figure}
See Fig. 1.  The total cost sending $p$ to $q$ is
\begin{equation}
 \label{eq:am420200228}
  C(p, q) = \int^1_0 \left|
   P^{-1}(z) - Q^{-1}(z)
  \right|^2 dz.
\end{equation}

We consider a regular statistical model
\begin{equation}
  S = \left\{ p(x, {\bm{\theta}}) \right\},
\end{equation}
parameterized by a vector parameter ${\bm{\theta}}$, where $p(x, {\bm{\theta}})$ is a probability density function of random variable $x \in {\bm{R}}^1$ with respect to the Lebesgue measure of ${\bm{R}}^1$.  Let
\begin{equation}
    D = \left\{ x_1, \cdots, x_n \right\}
\end{equation}
be $n$ independently observed data subject to $p(x, {\bm{\theta}})$.  We rearrange them in the increasing order,
\begin{equation}
    x_1 \le x_2 \le \cdots 
    \le x_n.
\end{equation}    
Then, $D$ is composed of order statistics.  We denote the empirical distribution by
\begin{equation}
    \hat{p}(x) = \frac 1n
    \sum \delta \left(x-x_i \right),
\end{equation}
where $\delta$ is the delta function.

The optimal transportation plan from $\hat{p}(x)$ to $p(x, {\bm{\theta}})$ is explicitly solved when $x$ is 1-dimensional, $x \in {\bm{R}}^1$.  The optimal plan is to transport a mass at $x$ to $x’$ defined by
\begin{equation}
    \hat{P}^{-1}(x) = P^{-1} 
    \left(x', {\bm{\theta}} \right),
\end{equation}
where $\hat{P}(x)$ and $P(x, {\bm{\theta}})$ are the cumulative distributions of $\hat{p}(x)$ and $p(x, {\bm{\theta}})$, respectively,
\begin{eqnarray}
 \hat{P}({\bm{x}}) &=& \int^x_{-\infty}
 \hat{p}(u) du, \\
 P(x, {\bm{\theta}}) &=& \int^x_{-\infty}
 p(u, {\bm{\theta}}) du,
\end{eqnarray}
and $\hat{P}^{-1}$, $P^{-1}$ are their inverse functions.
The total cost of transporting $\hat{p}(x)$ to $p(x, {\bm{\theta}})$ optimally is given by
\begin{equation}
    C({\bm{\theta}}) = \int^1_0
    \left|  
     \hat{P}^{-1}(z)-P^{-1}(z, {\bm{\theta}})
    \right|^2 dz.
\end{equation}

Let $z_1, \cdots, z_n$ be the points of equi-probability partition of $X$ for distribution $p(x, {\bm{\theta}})$ such that
\begin{equation}
 \label{eq:am1320200228}
    \int^{z_i}_{z_{i-1}}
    p(x, {\bm{\theta}}) dx =
    \frac 1n,
\end{equation}
where $z_0=-\infty$ and $z_n=\infty$.  In terms of the cumulative distribution, $z_i$ are written as
\begin{equation}
    P \left(z_i, {\bm{\theta}} \right)
    = \frac in
\end{equation}
and
\begin{equation}
    z_i = P^{-1} 
    \left( \frac in, {\bm{\theta}} \right).
\end{equation}
See Fig. 2.
\begin{figure}
    \centering
    \includegraphics[width=8cm]{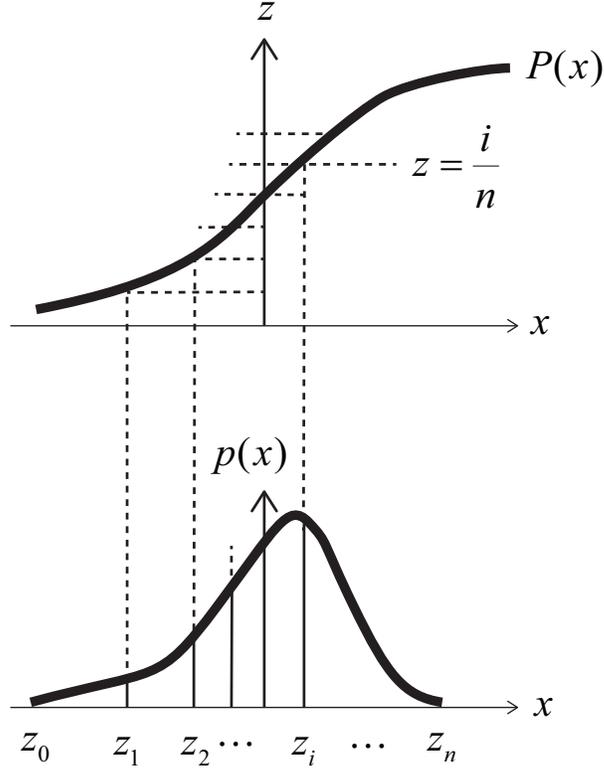}
    \caption{Equi-partition points $z_i$ of probability}
    \label{fig2}
\end{figure}
  
The optimal transportation cost is rewritten as
\begin{eqnarray}
  C({\bm{\theta}}) &=& \sum_i
  \int^{z_i}_{z_{i-1}}(x_i-z)^2
  p(z)dz \\
  &=& \frac 1n \sum x^2_i-2 \sum k_i
  ({\bm{\theta}})x_i + S({\bm{\theta}}),
\end{eqnarray}
where we use (\ref{eq:am1320200228}) and put
\begin{eqnarray}
  \label{eq:am1820200225}
  k_i({\bm{\theta}}) &=& 
  \int^{z_i}_{z_{i-1}} zp 
  (z, {\bm{\theta}})dz \\
  S({\bm{\theta}}) &=& 
  \sum \int^{z_i}_{z_{i-1}} z^2 p
  (z, {\bm{\theta}})dz =
  \int z^2 p(z, {\bm{\theta}})dz.
\end{eqnarray} 
By using the mean and variance of $p(x, {\bm{\theta}})$, 
\begin{eqnarray}
  \mu({\bm{\theta}}) &=& \int zp
  (z, {\bm{\theta}}) dz, \\
  \sigma^2({\bm{\theta}}) &=& \int z^2 p
  (z, {\bm{\theta}}) dz-\mu^2. 
\end{eqnarray}
We have
\begin{equation}
  S({\bm{\theta}}) = \mu^2 + \sigma^2.
\end{equation}

We define the $W$-estimator $\hat{\bm{\theta}}$ by the minimizer of $C({\bm{\theta}})$.  Differentiating $C({\bm{\theta}})$ with respect to ${\bm{\theta}}$ and putting it equal to 0, we have the estimating equation.

\begin{theorem}\upshape
The $W$-estimator $\hat{\bm{\theta}}$ satisfies
\begin{equation}
 \label{eq:am2320200225}
  \sum \frac{\partial}{\partial {\bm{\theta}}}
  k_i ({\bm{\theta}}) x_i =
  \frac 12 \frac{\partial}{\partial {\bm{\theta}}} S.
\end{equation}
\end{theorem}

It is interesting to see that the estimating equation is linear in $n$ observations $x_1, \cdots, x_n$ for any statistical model.  This is quite different from the maximum likelihood estimator or Bayes estimator.

We give a rough sketch that the estimator is asymptotically consistent, that is, it converges to the true ${\bm{\theta}}_0$ as $n$ tends to infinity. More detailed discussions are given for the location-scale model in the next section.  As $n$ tends to infinity, the order statistic $x_i$ converges to the $i$th partition point $z_i \left({\bm{\theta}}_0 \right)$, when the true parameter is ${\bm{\theta}}_0$.  From (\ref{eq:am1820200225}), we see that
\begin{equation}
    k_i = \frac 1n z_i({\bm{\theta}})
\end{equation}
as $n \rightarrow \infty$, so (\ref{eq:am2320200225}) is written as 
\begin{equation}
    \frac 1{2n} \frac{\partial}{\partial{\bm{\theta}}} \sum z_i ({\bm{\theta}}) z_i
    ({\bm{\theta}}_0)
    = \frac 12 S({\bm{\theta}}).
\end{equation}
We further remark that, as $n$ tends to infinity,
\begin{equation}
   \frac 1n \sum z^2_i = \int z^2 p
    (z, {\bm{\theta}}) dz =
    S({\bm{\theta}}).
\end{equation}
Therefore, ${\bm{\theta}}={\bm{\theta}}_0$ is the solution of (\ref{eq:am2320200225}) for $x_i=z_i \left({\bm{\theta}}_0\right)$, showing the consistency of the estimator.

\section{Location-scale model}

Let $f(x)$ be a standard probability density function, satisfying
\begin{eqnarray}
  \int f(x) dx &=& 1, \\
  \int xf(x) dx &=& 0, \\
  \int x^2 f(x) dx &=& 1,
\end{eqnarray}
that is, its mean is 0 and the variance is 1. The location-scale model $p(x, {\bm{\theta}})$ is written as
\begin{equation}
 \label{eq:am3120200225}
    p(x, {\bm{\theta}}) =
    \frac 1{\sigma} f
    \left( \frac{x-\mu}{\sigma} \right),
\end{equation}
where ${\bm{\theta}}=(\mu, \sigma)$ is the parameters to specify a distribution.

We define the equi-probability partition points $z_i$ for the standard $f(x)$ as
\begin{equation}
    z_i = F \left( \frac in \right),
\end{equation}
where $F$ is the cumulative distribution function
\begin{equation}
    F(x) = \int^x_{-\infty} f
    (u)du.
\end{equation}

We use the following transformation of the location and scale,
\begin{eqnarray}
  z &=& \frac{x-\mu}{\sigma}, \\
  x &=& \sigma z+ \mu.
\end{eqnarray}
The equi-probability partition points $\bar{x}_i$ of $p(x, {\bm{\theta}})$ is given by
\begin{equation}
    \bar{x}_i = \sigma z_i + \mu.
\end{equation}
The cost of the optimal transport from the empirical distribution $\hat{p}(x)$ to $p(x, \mu, \sigma)$ is then written as
\begin{eqnarray}
  C(\mu, \sigma) &=&  \sum
  \int^{\bar{x}_i}_{\bar{x}_{i-1}}
  \left(x_i-x \right)^2 p
  (x, \mu, \sigma) dx \nonumber\\
  \label{eq:am3720200225}
  &=&  
  \mu^2+\sigma^2 + \sum x^2_i -2 \sum
  x_i \int \left(\sigma z+ \mu \right)
  f(z)dz.
\end{eqnarray}
By differentiating (\ref{eq:am3720200225}), we obtain
\begin{eqnarray}
  \frac 12 \frac{\partial}{\partial \mu} C
   &=& \mu-\frac 1n \sum x_i, \\
   \frac 12 \frac{\partial}{\partial \sigma}C
   &=& \sigma-\sum k_i x_i,
\end{eqnarray}
where
\begin{equation}
 \label{eq:am4020200225}
    k_i = \int^{z_i}_{z_{i-1}} z f(z)dz,
\end{equation}
which does not depend on $\mu$ and $\sigma$ but depends only on the shape of $f$.  By putting the derivatives equal to 0, we obtain the following theorem.

\begin{theorem}\upshape
The $W$-estimator of a location-scale model is given by
\begin{eqnarray}
 \label{eq:am4120200225}
  \hat{\mu} &=& \frac 1n \sum x_i, \\
 \label{eq:am4220200225}
  \hat{\sigma} &=& \sum k_i x_i.
\end{eqnarray}
\end{theorem}

{\textbf{Remark}}\quad The $W$-estimator of the mean is the arithmetic average of observed data irrespective of the form of $f$.  The $W$-estimator of variance is also a linear function of observed data $x_1, \cdots, x_n$, but it depends on $f$, since $k_i$ depend on $f$.

The estimator $\hat{\mu}$ is consistent, asymptotically subject to the Gaussian distribution $N\left( \mu, \frac{\sigma^2}n \right)$.  We next show the asymptotic consistency of $\hat{\sigma}$ and its asymptotic variance.  Since the probability distribution of the order statistics $x_1, \cdots, x_n$ is explicitly given in literatures of statistics, it is, in principle, possible to calculate the variance, but we need complicated calculations.  So we here give a rough estimate based on speculative ideas.

\begin{theorem}\upshape
$\hat{\sigma}$ is asymptotically consistent with asymptotic variance
\begin{equation}
  V \left(\hat{\sigma}\right) =
  \frac{\sigma^2}n \int z^4 f(z)dz,
\end{equation}
where $V [\cdot]$ is the variance.
\end{theorem}

{\it{Sketch of proof.}} \quad
We evaluate $k_i$ when $n$ is large.  When $n$ is large, $z_{i-1}$ and $z_i$ are close and
\begin{equation}
    \Delta z_i = z_i-z_{i-1}
\end{equation}
is of order $1/n$.  More precisely, from
\begin{equation}
  \int^{z_i}_{z_{i-1}} f(z)dz = 
  \frac 1n,
\end{equation}
we have
\begin{equation}
  \Delta z_i f \left(z_i \right) =
  \frac 1n + O \left( \frac 1{n^2} \right).
\end{equation}
Hence, from (\ref{eq:am4020200225}), we have
\begin{equation}
    k_i = \frac 1n z_i + O 
    \left( \frac 1{n^2} \right).
\end{equation}
Thus, we have an asymptotic relation
\begin{equation}
 \label{eq:am4720200304}
 \hat{\sigma} =
 \frac{\sigma}n \sum z_i
 \hat{z}_i + \frac{\mu}n
 \sum \hat{z}_i,
\end{equation}
where
\begin{equation}
  \hat{z}_i =
  \frac{x_i-\mu}{\sigma}.
\end{equation}
We further use the following asymptotic relations
\begin{eqnarray}
  \frac 1n \sum z^2_i &\approx& \int
  z^2 p(z)dz = 1, \\
  \frac 1n \sum z_i &\approx& \int
  zp(z)dz = 0.
\end{eqnarray}
We finally have
\begin{equation}
  {\mathop{\lim}_{n \rightarrow \infty}}
  \hat{\sigma} = \sigma,
\end{equation}
showing that $\hat{\sigma}$ is asymptotically unbiased.

In order to evaluate the asymptotic variance, we use daring speculation.  To this end, we divide the $x$-axis into $n$ intervals $I_1 = \left[-\infty, z_1 \right], I_2 = \left[z_1, z_2 \right], \cdots, I_n = \left[z_{n-1}, \infty \right]$, the probability of each interval being equal to $1/n$.  When we select $n$ points from $f(x)$ independently, each observation $\hat{z}_i$ will fall into one interval randomly.  One interval may include multiple or no observations.  Let $s_i$ be a random variable to show the number of observations that fall in interval $I_i=\left[ z_{i-1}, z_i \right]$.  Then, each random variable $s_i$ is subject to Poisson distribution with mean and variance equal to 1.  They are independent except for the total constraint
\begin{equation}
  \sum s_i = n.
\end{equation}
The observed order statistic $\hat{z}_i$ will fall in interval $I_i = \left[z_{i-1}, z_i \right]$ most probably and takes value close to $z_i$.  It may fall in other nearby intervals.  

When $\hat{z}$, one of $\hat{z}$'s, falls in $I_i$, its value is written as
\begin{equation}
    \hat{z} = z_i-\varepsilon_i,
\end{equation}
where $\varepsilon_i$
\begin{equation}
  0 \le \varepsilon_i \le z_i-z_{i-1},
\end{equation}
is deviation within $I_i$.  It is a random variable of order $1/n$.  

Let us denote the interval $i'$ in which $\hat{z}_i$ falls.  Since $i$ and $i'$ are close, 
\begin{equation}
  \left| z_i-z_{i'} \right| =
  O \left( \frac 1n \right),
\end{equation}
with high probability, we can rewrite (\ref{eq:am4720200304}) as
\begin{equation}
  \hat{\sigma} = \frac{\sigma}n \sum
  s_{i'} z_{i'} \hat{z}_{i'} +
  O \left( \frac 1n \right)
\end{equation}
by neglecting high-order terms, where summation with respect to $i$ is replaced by summation with respect to the intervals $I_{i'}$ with weight $s_{i'}$.
When $s_i=0$, interval $I_i$ includes no observation.  When $s_i>1$, $I_i$ includes multiple observations.

We calculate the variance of (\ref{eq:am4220200225}) as
\begin{equation}
  V \left[ \hat{\sigma} \right] = V
  \left[ \frac{\sigma}n \sum_i s_i z^2_i \right]
  + O \left( \frac 1{n^2} \right).
\end{equation}
We further note that $s_i$ are asymptotically independent.  Hence, we have
\begin{eqnarray}
  V \left[\hat{\sigma} \right]
  &\approx& \frac{\sigma^2}{n^2} \sum V
  \left[s_i \right] z^4_i \\
  &\approx& \frac{\sigma^2}n \int z^4
  f(z)dz,
\end{eqnarray}
proving the theorem.

It is easy to see from (\ref{eq:am4120200225}) and (\ref{eq:am4220200225}) that $\hat{\mu}$ and $\hat{\sigma}$ are asymptotically non-correlated, since $x_i$'s are independent.

When $f$ is Gaussian
\begin{equation}
  f(z) = \frac 1{\sqrt{2 \pi}} \exp
  \left\{ -\frac 12 z^2 \right\},
\end{equation}
the asymptotic variance is
\begin{equation}
  V \left[\hat{\sigma}\right] =
  \frac 3n \sigma^2.
\end{equation}
Hence, it is consistent but not efficient.

When $f$ is uniform,
\begin{equation}
  f(z) = \left\{
  \begin{array}{cc}
   \frac 1{2 \sqrt{3}}, & 
   |z| \le \sqrt{3}, \\
    0, &  \mbox{otherwise}, 
  \end{array}
  \right.
\end{equation}
the asymptotic variance is
\begin{equation}
  V \left[\hat{\sigma}\right] =
  \frac 9{5n} \sigma^2.
\end{equation}
However, the Fisher information divergence to infinity for the uniform distribution and the maximum likelihood estimator $\hat{\sigma}$
converges to 0 exponentially fast.

In general, the $W$-estimator is not sensitive to changes of the waveform $f$, whereas the maximum likelihood estimator is sensitive.

\section{Riemannian structure of $W$-divergence}

Consider the manifold $M = \{p(x) \}$ of probability distributions which are absolutely continuous with respect to the Lebesgue measure and have finite second moments.  It is known that $M$ has Riemannian structure due to the Wasserstein distance or the cost function. For two distributions $p(x)$ and $q(x)$, their optimal transportation cost, that is, the divergence between them, is given by (\ref{eq:am420200228}).

We calculate the optimal transportation cost between two nearby distributions $p(x)$ and $p(x)+ \delta p(x)$, where $\delta p(x)$ is infinitesimally small.  We have
\begin{equation}
  \left(P + \delta P \right)^{-1} (z)
  = P^{-1}(z)- 
  \frac{\delta P \left\{ x(z) \right\}}{P' \left\{ x(z) \right\}},
\end{equation}
where
\begin{equation}
  x(z) = P^{-1}(z).
\end{equation}
This equation is derived from
\begin{equation}
  \frac d{dz} F^{-1}(z) =
  \frac 1{f' \left\{ x(z) \right\}},
\end{equation}
which we have from the differentiation of the identity
\begin{equation}
  F^{-1} \left\{ F(x) \right\} = x.
\end{equation}
We thus have
\begin{equation}
 \label{eq:am6920200225}
  C \left(p, p+\delta p \right) =
  \int^{\infty}_{-\infty}
  \frac 1{p(x)} \left(
   \int^{x}_{-\infty} \delta p(y)dy
  \right)^2 dx
\end{equation}
which is a quadratic form of $\delta p(x)$.  This gives a Riemannian metric to $M$.

The location-scale model $S$ is a finite-dimensional submanifold embedded in $M$.  We have for the location-scale model (\ref{eq:am3120200225}),
\begin{equation}
  \delta p(y) = \frac{\partial}{\partial \mu}
  p(y, {\bm{\theta}}) d \mu +
  \frac{\partial}{\partial \sigma}
  p(y, {\bm{\theta}}) d \sigma.
\end{equation}
The Riemannian metric tensor $G=\left(g_{ij}\right)$ is derived from
\begin{equation}
  C (p, p+\delta p) = \sum g_{ij}
  ({\bm{\theta}}) d \theta_i d \theta_j.
\end{equation}
See also \citet*{LZ2019}.

\begin{theorem}\upshape
The location-scale model is a Euclidean space, irrespective of $f$,
\begin{equation}
  g_{ij} = \delta_{ij}.
\end{equation}
\end{theorem}

\begin{proof}\upshape
We need to calculate (\ref{eq:am6920200225}).  Technical details are given in Appendix.
\end{proof}
It is surprising that $G=\left(g_{ij}\right)$ is the identity matrix for the location-scale model, so that $S$ is a Euclidean space.  See also \citealt*{LZ2019}.  It is flat by itself, but $S$ is a curved submanifold in $M$ \citep{Takatsu2011}, like a cylinder embedded in ${\bm{R}}^2$.

When $n$ is large, the cost decreases in the order of $1/n$.  The $W$-estimator is the projection of $\hat{p}(x)$ to $S$ in the tangent space of $M$.  Let $\hat{\theta}'$ be another consistent estimator. Then, we have the Pythagorean relation
\begin{equation}
  C \left(\hat{p}, p_{\hat{\theta}'} \right)
  = C \left(\hat{p}, p_{\hat{\theta}} \right)
  + C \left(p_{\hat{\theta}}, p_{\hat{\theta}'} \right)
\end{equation}
and the difference of the cost between the two estimators is 
\begin{equation}
 C \left( p_{\hat{\theta}}, p_{\hat{\theta}'} \right)
 = \left|{\bm{\theta}}'- \hat{\bm{\theta}}' \right|^2.
\end{equation}

\citet*{LZ2019} studies the properties of the W estimator given by the $W$ score function. They give the $W$-efficiency and $W$ Cramer-Rao inequality.  However, their $W$-estimator does not minimize the transportation cost.  It is interesting to study the relation between the two $W$-estimators.

\section{Conclusions}

We studied the behaviors of the $W$-estimator minimizing the transportation cost from the observed empirical distribution to the underlying statistical model on ${\bm{R}}^1$.  It is a consistent estimator having a simple form of the estimating equation.  We focused on the location-scale model and showed that the estimator can be represented by a simple linear form of observations.  Its asymptotic variance was calculated.  Although its error variance is worse than the maximum likelihood estimator, it is simple, and further it is the estimator that minimizes the transportation cost from the observed sample to the model.

We need to study further its merits and demerits.  We hope to find good applications to computer vision and AI.  It is an interesting problem to compare the $W$-estimator of \citet*{LZ2019} which uses the $W$ score function with the minimum cost $W$-estimator.

\section*{Appendix: The Riemannian metric of the location scale model}

We have
\begin{equation}
  \delta p(x, {\bm{\theta}}) = 
  -\frac 1{\sigma^2} f' \left( 
  \frac{x-\mu}{\sigma} \right) d \mu
  -\frac 1{\sigma^3} \left\{
   \sigma f \left( \frac{x-\mu}{\sigma} \right)
   + (x-\mu) f' \left( \frac{x-\mu}{\sigma} \right)
  \right\} d \sigma.
\end{equation}
By integration, we have
\begin{equation}
  \int^{x}_{-\infty} \delta p
  (y, {\bm{\theta}}) dy =
  -p(x, {\bm{\theta}}) d \mu -
  (x-\mu) p(x, {\bm{\theta}}) d \sigma.
\end{equation}
Hence, we have
\begin{equation}
  C({\bm{\theta}}, {\bm{\theta}}+ d{\bm{\theta}})
  = d \mu^2 + d \sigma^2.
\end{equation}


\begin{thebibliography}{}
%
\bibitem[Amari(2016)]{Amari2016} Amari, S., 
Information Geometry and Its Applications. Springer (2016).
%
\bibitem[Amari et al.(2018) Amari, Karakida and Oizumi]{AKO2018} Amari, S., Karakida, R., Oizumi, M., Information geometry connecting Wasserstein distance and Kullback-Leibler divergence via the entropy-relaxed transportation problem. Information Geometry, 1, 13--37, (2018).
%
\bibitem[Amari et al.(2019) Amari, Karakida, Oizumi and Cuturi]{AKOC2019} Amari, S., Karakida, R., Oizumi, M., Cuturi, M., Information geometry for regularized optimal transport and barycenters of patterns. Neural Computation, 31, 827--848, (2019).
%
\bibitem[Arjovsky et al.(2017) Arjovsky, Chintala and Bottou]{ACB2917} Arjovsky, M., Chintala, S., Bottou, L., Wasserstein GAN. arXiv:1701.07875, (2017).

%
\bibitem[Fronger et al.(2015) Fronger, Zhang, Mobahi, Araya-Polo and Poggio]{FZMAP2015} Fronger, C. Zhang, C., Mobahi, H., Araya-Polo, M., Poggio, T., Learning with a Wasserstein loss. NIPS, 28, (2015).
%
\bibitem[Kurose et al.(2019) Kurose, Yoshizawa and Amari]{KYA2019} Kurose, T., Yoshizawa, S. and Amari, S., Optimal transportation plan with generalized entropy regularization. submitted, (2019).
%
\bibitem[Li et al.(2018) Li and Mont\'{u}far]{LM2018} Li, W., Mont\'{u}far, G., Ricci curvature for parametric statistics via optimal transport. 	arXiv:1807.07095 (2018).
%
\bibitem[Li et al.(2019) Li and Zhao]{LZ2019} Li, W., Zhao, J., Wasserstein information matrix. memo (2019).  
%
\bibitem[Montavon et al.(2015) Montavon, Muller and Cuturi]{MMC2015} Montavon, G., Muller, K., Cuturi, M., Wasserstein training for Boltzmann machine. aeXiv:1507.01972v1, (2015).
%
\bibitem[Peyr\'{e} et al.(2019) Peyr\'{e} and Cuturi]{PC2019} Peyr\'{e}, G., Cuturi, M., Computational optimal transport (2019).
%
\bibitem[Takatsu(2011)]{Takatsu2011} Takatsu, A., Wasserstein geometry of Gaussian measures. Osaka J. Math., 48, 1005--1026, (2011).
%
\bibitem[Villani(2003)]{Villani2003} Villani, C., Topics in Optimal Transportation. American Mathematical Society, (2003).
%
\bibitem[Villani(2009)]{Villani2009} Villani, C., Optimal Transport, Old and New. Springer, (2009).
%
\bibitem[Wang et al.(2019) Wang and Li]{WL2019} Wang, Y., Li, W., Information Newton’s flow: Second-order optimization method in probability space. arXiv, (2019).
\end{thebibliography}
\end{document}